\documentclass[11pt]{amsart}
\usepackage{amsthm, amsmath,amssymb, psfrag}
\usepackage[curve,matrix,arrow]{xy}
\usepackage{mathrsfs}
\usepackage{isotope}
\renewcommand{\isotopestyle}{\it}
\usepackage{graphics}
\DeclareMathOperator{\Ho}{Ho}
\DeclareMathOperator{\Id}{Id}

\DeclareMathOperator{\colim}{colim}

\DeclareMathOperator{\ord}{-ord}

\newcommand{\mF}{{\mathbb F}}
\newcommand{\mQ}{{\mathbb Q}}
\newcommand{\mR}{{\mathbb R}}
\newcommand{\mZ}{{\mathbb Z}}
\newcommand{\C}{{\mathcal C}}
\newcommand{\Ac}{{\mathcal A}}
\newcommand{\Dc}{{\mathcal D}}

\newcommand{\Rc}{{\mathcal R}}
\newcommand{\Sc}{{\mathcal S}}
\newcommand{\Tc}{{\mathcal T}}
\newcommand{\iso}{\cong}
\newcommand{\xra}{\xrightarrow}

\newcommand{\sm}{\wedge}
\newcommand{\SW}{\mathcal S \mathcal W}
\newcommand{\tensor}{\otimes}

\newcommand{\td}[1]{\langle #1\rangle}
\renewcommand{\to}{\longrightarrow}

\newtheorem{theorem}[equation]{Theorem}

\newtheorem{prop}[equation]{Proposition}

\theoremstyle{definition}
\newtheorem{defn}[equation]{Definition}

\newtheorem{rk}[equation]{Remark}
\newtheorem{problem}[equation]{Open problem}

\begin{document}

\title{\mbox{Algebraic versus topological triangulated categories}}

\date{\today; 2000 AMS Math.\ Subj.\ Class.: 18E30, 55P42}
\author{Stefan Schwede}
\address{Mathematisches Institut, Universit\"at Bonn, Germany}
\email{schwede@math.uni-bonn.de}
\maketitle

These are extended and updated notes of a talk,
the first version of which I gave at the 
{\em Workshop on Triangulated Categories} at the University of Leeds,
August 13-19, 2006. These notes are mostly expository and do not contain
all proofs; I intend to publish the remaining details elsewhere.

The most commonly known triangulated categories arise from chain 
complexes in an abelian category by passing to chain homotopy classes
or inverting quasi-isomorphisms. Such examples are called `algebraic' 
because they originate from abelian (or at least additive) categories.
Stable homotopy theory produces examples of triangulated categories
by quite different means, and in this context the source categories
are usually very `non-additive' before passing to homotopy classes of 
morphisms. Because of their origin I refer to these examples as 
`topological triangulated categories'. 

In this note I want to explain some systematic differences between these 
two kinds of triangulated categories. There are certain properties 
-- defined entirely in terms of the triangulated structure --
which hold in all algebraic examples, but which fail in some topological 
ones. These differences are all torsion phenomena, and rationally there is 
no difference between algebraic and topological triangulated categories.

\bigskip

A triangulated category is {\em algebraic} in the sense of 
Keller~\cite[3.6]{keller-differential graded} if it is triangle equivalent
to the stable category of a Frobenius category, i.e., an exact category
with enough injectives and enough projectives in which injectives
and projectives coincide. Examples include all triangulated categories
which one should reasonably think of as `algebraic': various homotopy
categories and derived categories of rings, schemes and abelian categories;
stable module categories of Frobenius rings; derived categories of
modules over differential graded algebras and differential graded categories.
By a theorem of Porta~\cite[Thm.~1.2]{porta}, 
every algebraic triangulated category which is {\em well generated} 
(a mild restriction on its `size', see~\cite[Def.~8.1.6 and 8.1.7]{neeman}) 
is equivalent to a localization of the derived category $\Dc(\Ac)$
of a small differential graded category $\Ac$.

For an object $X$ of a triangulated category $\Tc$ and a 
natural number $n$ we write $n\cdot\Id_X$ or simply $n\cdot X$ for the $n$-fold
multiple of the identity morphism in the group $[X,X]$
of endomorphisms in $\Tc$.
We let $X/n$ denote any cone of $n\cdot\Id_X$,
i.e., any object which is part of a distinguished triangle
$$ X\ \xra{\ n\cdot\, }\ X\ \xra{\quad}\ X/n\ \xra{\quad}\ X[1] \ . $$
A short diagram chase shows that the group $[X/n,X/n]$ is 
always annihilated by $n^2$; 
in algebraic triangulated categories, more is true:

\begin{prop}\label{prop-mod n in algebraic}
  If $\Tc$ is algebraic, then $n\cdot X/n=0$.  
\end{prop}
\begin{proof}
  We exploit that algebraic triangulated categories are
tensored over $\Dc^b(\mZ)$, the bounded derived category of finitely generated
abelian groups.
This means that there is a biexact functor
$$ \tensor^L \ : \ \Dc^b(\mZ) \times \Tc \ \to \ \Tc$$
which is associative and unital up to coherent isomorphism
with respect to the derived tensor product in $\Dc^b(\mZ)$.
In  $\Dc^b(\mZ)$ we have a distinguished triangle
$$ \mZ\ \xra{\ n\cdot\ }\ \mZ\ \xra{\quad}\ \mZ/n\ \xra{\quad}\ \mZ[1]$$
which becomes a distinguished triangle in $\Tc$ after tensoring
with $X$. So $\mZ/n\tensor^L X$ is isomorphic to $X/n$.
Since the tensor product is additive in each variable and
since $n\cdot\mZ/n=0$ in $\Dc^b(\mZ)$, 
we conclude that $n\cdot(\mZ/n\tensor^L X)=0$.
\end{proof}

In contrast to Proposition~\ref{prop-mod n in algebraic}, 
in a general triangulated category we
can have $n\cdot X/n\ne 0$ for suitable choices of $X$ and $n$.
An example arises in the {\em Spanier-Whitehead category}
which we now review. The Spanier-Whitehead category is made from
homotopy classes of continuous maps between certain kinds
of topological spaces, and it is the prime example
of a {\em topological} triangulated category (to be defined below)
which is not algebraic. The Spanier-Whitehead category was originally
introduced (without the formal desuspensions) in~\cite{spanier-whitehead}.

We recall that the {\em reduced suspension} $\Sigma X$ of a space $X$ 
with basepoint $x_0$ is given by\vspace{.8cm}
\[
\Sigma X = \frac{X\times[0, 1]}{X\times\{0, 1\}\cup \{x_0\}\times[0, 1]}\ .
\hspace*{5cm}\raisebox{-7ex}[1.5ex][2.0ex]{\scalebox{.60}{\includegraphics{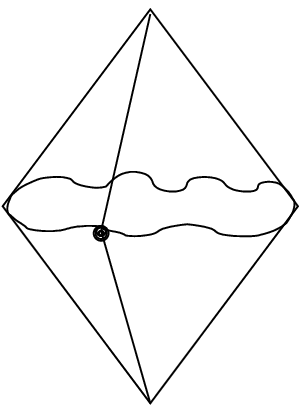}}}
\]
\vspace{.8cm}\\
For example we have $\Sigma S^{n-1}\iso S^n$, i.e., the suspension of 
a sphere is homeomorphic to a sphere of the next dimension.

\begin{defn}
The {\em Spanier-Whitehead category}, denoted $\SW$,
has as objects the pairs $(X, n)$ where $X$ is a finite
CW-complex equipped with a distinguished basepoint 
and $n\in\mZ$. Morphisms are given by
\[  \SW ((X, n), (Y, m)) \  =\ 
\colim_{k\rightarrow\infty} [\Sigma^{k+n} X, \Sigma^{k+m} Y] \]
where square brackets $[-,-]$ denote pointed homotopy classes of
continuous, basepoint preserving maps.
The colimit is formed by iterated suspensions; by Freudenthal's 
suspension theorem, it is actually attained at a finite stage. 
Composition in $\SW$ is defined by composition of representatives, suitably
suspended so that composition is possible.
\end{defn}

It is convenient to identify a finite pointed CW-complex $X$ with 
the object $(X,0)$ of the Spanier-Whitehead category.
Then two CW-complexes become isomorphic in $\SW$ if and only 
if they become homotopy equivalent after a finite number of suspensions.

The Spanier-Whitehead category is triangulated in the following way:
\begin{itemize}
\item The shift functor is given by $(X,n)[1]=(X,n+1)$.
Tautologically, the identity map of $\Sigma^{n+1} X$ is an isomorphism
between $(\Sigma X, n)$ and $(X, n+1)$ in $\SW$,
so suspension is in fact isomorphic to the shift and is invertible in $\SW$.
\item The Spanier-Whitehead category is additive:
for pointed spaces $X$ and $Y$, the
set $[\Sigma X, Y]$ of pointed homotopy classes from a suspension
has a natural group structure as follows. 
The product of $f, g : \Sigma X\to Y$ is represented by the composite
\begin{center}
\scalebox{.60}{\includegraphics{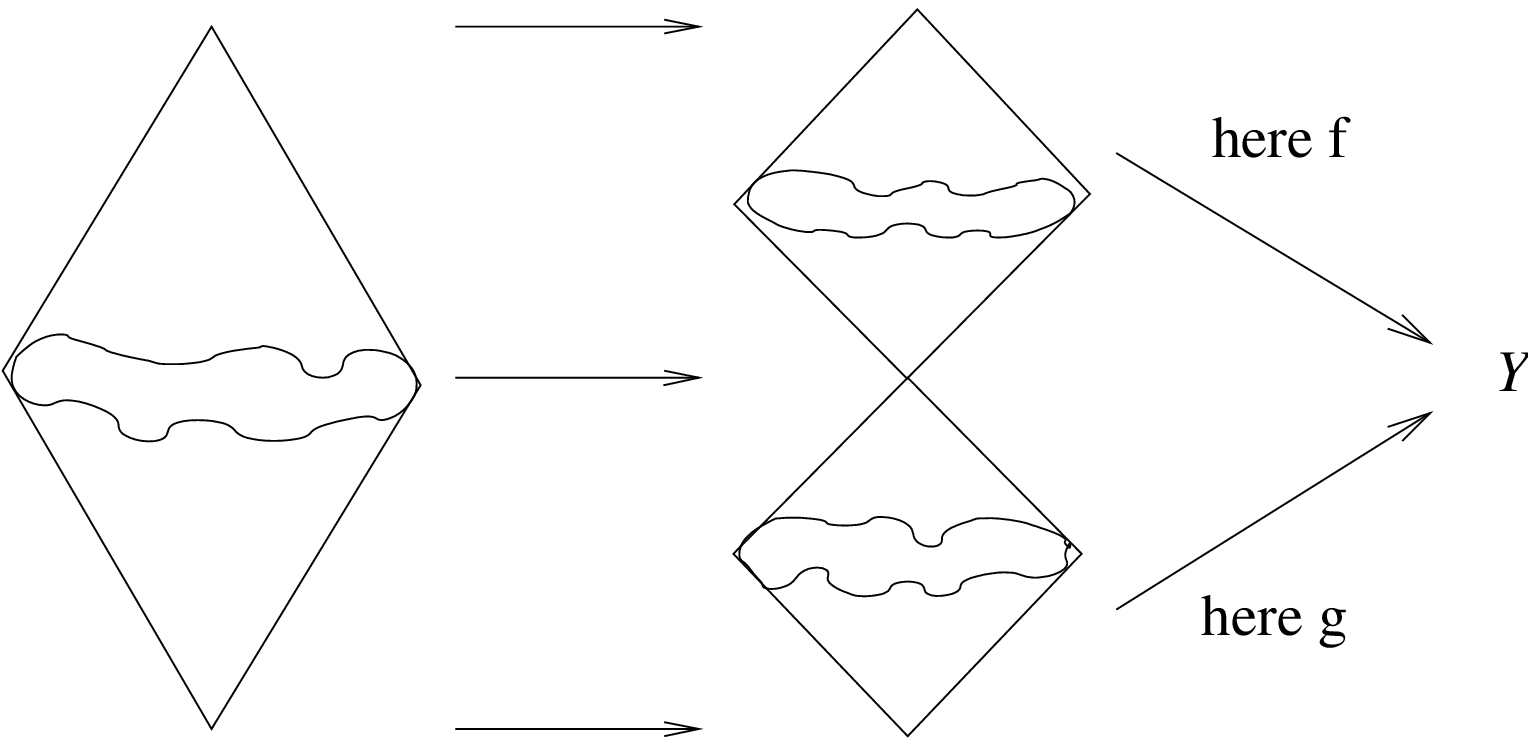}}
\end{center}
\[\xymatrix@C=30mm{
& \Sigma X \ \ar[r]^-{\text{pinch}} &
\ \Sigma X\vee\Sigma X\ \ar[r]^{\  f\vee g\ } &\  Y &}
\]
For $X=S^0$ we have $\Sigma X\iso S^1$ and this reduces to the group structure 
on the fundamental group $[S^1, Y] = \pi_1 (Y, y_0)$. 
On a double suspension as source object, this group structure is abelian;
an example of this is that for $n\geq 2$ the higher homotopy groups
$\pi_n (Y, y_0) = [S^n, Y]$ are abelian.
In the Spanier-Whitehead category, every object is a double suspension, 
so the homomorphism sets in $\SW$ are naturally abelian groups.
\item Mapping cone sequences give distinguished triangles:
the {\em mapping cone} of a
pointed map $f : X\to Y$ is the space
\[
C(f) = \frac{X\times[0, 1]\cup_{X\times \{1\}} Y}
{X\times \{0\}\cup \{x_0\}\times [0, 1]}\qquad\qquad\raisebox{-6ex}[9ex][6.0ex]
{\scalebox{.40}{\includegraphics{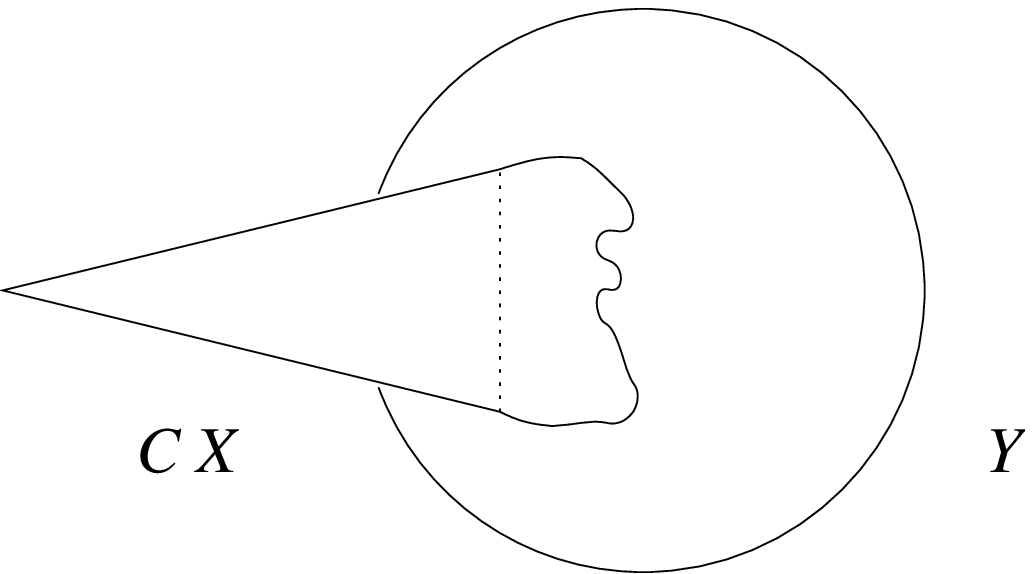}}}
\]
There is an inclusion $i:Y\to C(f)$ and a projection
$p:C(f)\to\Sigma X$ (which collapses $Y$ to a point). 
In the sequence of pointed maps
\[
X\ \xra{\  f\ }\ Y\ \xra{\ i\ }\ C(f)\ \xra{\ p\ }\ \Sigma X
\]
the composite of any two is null-homotopic, and these sequences
and their \mbox{(de-)suspensions} give the distinguished triangles in $\SW$.
\end{itemize}

A verification of the axioms of a triangulated category 
(except for the octahedral axiom), and more details
on the Spanier-Whitehead category can be found in Ch.~1, \S 2
of Margolis book~\cite{margolis}. (Margolis does not impose any finiteness
restriction on the objects of the Spanier-Whitehead category;
the category $\SW$ is denoted {\bf SW}$_{\text{f}}$ in~\cite{margolis}.)

The group
\[
\SW (S^n, X)\ =\ \colim_{k\rightarrow\infty}\  [S^{k+n}, \Sigma^k X] \ = \
\colim_k \,\pi_{k+n}(\Sigma^k X,*)
\]
is denoted by $\pi_n^s X$ and called the $n$-th {\em stable homotopy group}
of~$X$. For $X=S^0$, we abbreviate $\pi_n^sS^0$ to $\pi_n^s$ and
speak about the {\em $n$-th stable homotopy group of spheres}, 
also called the {\em $n$-th stable stem}.
For example, $\pi_0^s=\mZ$, generated by the identity map, and
$\pi_1^s=\mZ /2$ generated by the class of the Hopf map $\eta:S^3\to S^2$.
The stable stems are easy to define, but notoriously hard to compute;
for example, there is no finite CW-complex $X$ which is non-trivial in $\SW$
for which all stable homotopy groups are known! Much
machinery of algebraic topology has been developed to calculate such groups
and understand their structure, but no one expects to ever get explicit
formulae for all stable homotopy groups of spheres.

We also need the symmetric monoidal {\em smash product} 
on the Spanier-Whitehead category. This arises from the geometric 
smash product $X\sm Y$ of two pointed spaces $X$ and~$Y$
defined as
$X\sm Y=(X\times Y)/(X\times\{y_0\}\cup\{x_0\}\times Y)$.
Suspension is an example of the smash product,
i.e., $\Sigma X$ is naturally homeomorphic to $X\sm S^1$.
The smash product can be extended to the Spanier-Whitehead category by
$(X,n)\sm(Y,m)=(X\sm Y,n+m)$, and then it becomes biexact, 
i.e., an exact functor of triangulated categories in each variable.

We denote by $S=(S^0,0)$ the unit object  of the smash product in $\SW$
and refer to it as the {\em sphere spectrum}. We use this terminology
because the Spanier-Whitehead category can be identified with the
compact objects in a larger triangulated category,
the {\em stable homotopy category}, whose objects are called spectra.
However, the differences between algebraic
and topological triangulated categories already show up in the
Spanier-Whitehead category, which is easier to define than the
full stable homotopy category.

For $n\geq 2$, the {\em mod-$n$ Moore spectrum} is defined 
as a cone of multiplication by $n$ on the sphere spectrum,
i.e., it is part of a distinguished triangle
\begin{equation}\label{eg-triangle defining Moore} 
S \ \xra{\ n\cdot\, }\ S \ \xra{\quad} \ S/n \ \xra{\quad} \ S[1] \ . 
\end{equation}
More concretely we can define $S/n=(S^1\cup_n D^2,-1)$, 
the formal desuspension of a two-dimensional mod-$n$ Moore space, 
i.e., the space obtained from the
circle $S^1$ by attaching a 2-disc along the degree~$n$ map
$S^1\to S^1, z\mapsto z^n$. For example, $S/2=(\mR P^2,-1)$,
the formal desuspension of 2-dimensional real projective space.
If $n$ and $m$ are coprime, then $S/(nm)$ is isomorphic
in the Spanier-Whitehead category to the smash product $S/n\sm S/m$.
So we often concentrate on Moore spectra for primes or prime powers,
since the other Moore spectra are smash products of these.

Proposition~\ref{prop-mod n in algebraic}  and the following 
proposition show that the Spanier-Whitehead category is not algebraic.

\begin{prop}\label{prop-mod 2 Moore}
The morphism $2\cdot S/2$ is nonzero in $\SW$.
\end{prop}

The standard tool for proving Proposition~\ref{prop-mod 2 Moore} and
its generalizations below are cohomology operations, which we quickly review.
The {\em mod-$p$ cohomology} of an object $(X,n)$ in $\SW$ defined by 
$$ H^k((X,n),\mF_p) \ =\ \tilde H^{k-n}(X,\mF_p)$$
where the right hand side is reduced singular cohomology
with coefficients in the field~$\mF_p$.
These cohomology groups have a natural action (essentially by definition)
by the mod-$p$ {\em Steenrod algebra}, i.e., the algebra of stable,
natural, graded mod-$p$ cohomology operations. The Steenrod algebra
is generated for $p=2$ by operations $Sq^i$ of degree~$i$ for $i\geq 1$
and for odd $p$ by the Bockstein operation $\beta$ of degree~1
and operations $P^i$ of degree $i(2p-2)$ for $i\geq 1$. For $p=2$, the operation
$Sq^1$ equals the Bockstein operation. These operations satisfy
the {\em Adem relations}, which we do not reproduce here.

The mod-$n$ Moore spectrum is characterized up to 
isomorphism in the Spanier-Whitehead category by the property
that its integral spectrum homology is
concentrated in dimension zero where it is isomorphic to $\mZ/n$.
For a prime~$p$ the mod-$p$ cohomology of $S/p$ is
one-dimensional in dimensions~0 and~1, and trivial otherwise,
and the Bockstein operation is non-trivial from dimension~0 to
dimension~1.

\begin{proof}[Proof of Proposition~\ref{prop-mod 2 Moore}]
This proposition is a classical fact, which should probably be credited
to Steenrod since the standard proof uses mod-2 Steenrod operations.
We argue by contradiction and suppose that $2\cdot S/2=0$.
If we smash the defining triangle~\eqref{eg-triangle defining Moore} 
with another copy of the mod-2 Moore spectrum and use that $S$
is the unit of the smash product, we obtain a distinguished 
triangle
$$ S/2\ \xra{\ 2\cdot\, }\ S/2\ \xra{\quad}\ S/2\sm S/2
\ \xra{\quad}\ S/2[1] \ . $$
Under our assumption the first map is trivial so the smash product
$S/2\sm S/2$ splits in $\SW$ as a sum of $S/2$ and $S/2[1]$.
Thus as a module over the Steenrod algebra,
the mod-2 cohomology of the smash product
$S/2\sm S/2$ decomposes into a sum of two non-trivial summands.

On the other hand, there is a K\"unneth isomorphism for the
 mod-2 cohomology of a smash product 
$$ H^n(X\sm Y,\mF_2) \ \iso \ \bigoplus_{p+q=n}\,  
 H^p(X,\mF_2)\tensor H^q(Y,\mF_2) \ . $$
This is an isomorphism as modules over the Steenrod algebra
with action on the right hand side given by the {\em Cartan formula}
$$ Sq^m(x\tensor y) \ = \ \sum_{i=0}^m \, Sq^ix\tensor Sq^{m-i}y \ .$$
In the mod-2 cohomology of the mod-2 Moore spectrum the
Bockstein operation\linebreak
\mbox{$Sq^1:H^0(S/2,\mF_2)\to H^1(S/2,\mF_2)$} is non-zero.
The Cartan formula shows that the operation $Sq^2$ is then non-trivial
on the tensor product of two copies of the generator of $H^0(S/2,\mF_2)$.
Thus the mod-2 cohomology of $S/2\sm S/2$ is a 4-dimensional,
indecomposable module over the Steenrod algebra. We have reached
a contradiction, which means that we must have  $2\cdot S/2\ne 0$
in $\SW$.
\end{proof}

In topological triangulated categories, the phenomenon that we
can have $n\cdot X/n\ne 0$ is entirely 2-local. To explain this, 
I have to be more precise about what I mean by
a topological triangulated category. 
A {\em model category} (in the sense of Quillen~\cite{Q}) is 
an axiomatic framework for homotopy theoretic constructions.
Among other things, a model category structure allows one 
to define mapping cones and suspensions and talk 
about homotopies between morphisms. A model category is called
{\em stable} if it has a zero object and the suspension functor 
is a self-equivalence of its homotopy category. 
The homotopy category of a stable model category is naturally 
a triangulated category (cf.~\cite[7.1.6]{hovey-book});
the proof is essentially the same as for the Spanier-Whitehead category.
By definition the suspension
functor is a self-equivalence, and it defines 
the shift functor. Since every object is a two-fold suspension,
hence an abelian co-group object, the homotopy category of
a stable model category is additive.
The distinguished triangles are defined by mapping cone sequences.

For us a {\em topological triangulated category}
is any triangulated category which is equivalent
to a full triangulated subcategory of the homotopy category
of a stable model category. An important example which was already mentioned
above is the 
{\em stable homotopy category} of algebraic topology which 
was first introduced by Boardman (unpublished;
accounts of Boardman's construction appear in~\cite{vogt-aarhus}
and~\cite[Part~III]{adams-generalized}). 
There is an abundance of models for the
stable homotopy category, see for 
example~\cite{BF, EKMM, HSS, MMSS, lydakis-simplicial}.
The Spanier-Whitehead category $\SW$ is equivalent to the full subcategory
of compact objects in the stable homotopy category,
so it is a topological triangulated category in our sense.
Further examples of topological triangulated categories 
are `derived' (i.e., homotopy) categories of structured ring spectra, 
equivariant and motivic stable homotopy categories,
sheaves of spectra on a Grothendieck site or (Bousfield-) localizations
of all these, see~\cite[Sec.~2.3]{ss-modules} for more details.
The theorem of Porta mentioned above has an analogue in this context:
Heider essentially shows in~\cite[Thm.~4.7]{heider} 
that every topological triangulated category 
which is well generated is equivalent to a localization 
of the homotopy category $\Ho(\Rc\text{-mod})$
of a small spectral category $\Rc$.

Algebraic triangulated categories are typically also topological
(the converse is not generally true, and that is the point of these notes).
For algebraic triangulated categories which are derived categories
of abelian categories, this follows whenever there is a model
structure on the category of chain complexes with quasi-isomorphisms
as weak equivalences (see for example~\cite{hovey-sheaves} 
or Section~2.4 of~\cite{ss-modules}
for more details and references). Similarly, for modules over
a Frobenius ring, there is a stable model structure with stable
equivalences as the weak equivalences, see~\cite[Thm.~2.2.12]{hovey-book}.
More generally, any well-generated algebraic triangulated category 
is equivalent to a localization of the derived category $\Dc(\Ac)$
of small differential graded category $\Ac$~\cite[Thm.~1.2]{porta}. 
The localization can be realized as a Bousfield localization of the ordinary
(i.e., `projective') model structure on modules over $\Ac$;
hence $\Dc(\Ac)$ and its localizations are topological.

Examples of triangulated categories which are 
neither algebraic nor topological were recently constructed by
Muro, Strickland and the author~\cite{nomodel}. The simplest one
is the category $\mathcal F(\mZ/4)$ of finitely generated free modules over 
the ring $\mZ/4$. The category $\mathcal F(\mZ/4)$ has a unique triangulation with
the identity shift functor and such that the triangle
$$ \mZ/4\ \xra{\ 2\ }\ \mZ/4\ \xra{\ 2\ }\ \mZ/4\ \xra{\ 2\ }\ \mZ/4$$
is exact. For the proof of this and an argument why  $\mathcal F(\mZ/4)$
is not topological I refer to~\cite{nomodel}. At present,
I do not know any `exotic' (i.e., non-topological and non-algebraic)
triangulated category in which~2 is invertible.

Unlike algebraic triangulated categories, we can not usually expect 
that a topological triangulated category 
can be tensored over $\Dc^b(\mZ)$, the bounded derived category of
finitely generated abelian groups. 
The appropriate replacement for $\Dc^b(\mZ)$ is the Spanier-Whitehead category:
for every topological triangulated category $\Tc$, 
there is a biexact pairing
$$ \sm \ : \ \SW \times \Tc \ \to \ \Tc$$
which is associative and unital up to coherent natural isomorphism
with respect to the smash product in the Spanier-Whitehead category.
The shift functor in $\Tc$ is isomorphic to smashing with the circle $S^1$,
view as an object of $\SW$.

We use the above pairing as a black box, but I'll briefly indicate how
it can be constructed.
One way is to start from the action of the homotopy
category of pointed simplicial sets on the homotopy category of
a pointed model category $\C$,
$$ \sm \ : \ \text{Ho{\bf SSet}}_* \times \text{Ho }\C \ 
\to \ \text{Ho }\C \ , $$
which uses a technique called `framings', see~\cite[Thm.~5.7.3]{hovey-book}. 
The geometric realization of a finite simplicial set is a finite CW-complex, 
and up to homotopy equivalence, every finite CW-complex arises in this way.
Since moreover suspension is invertible in Ho~$\C$,
this extends to a well-defined pairing on the Spanier-Whitehead category.

There is an alternative construction for {\em spectral model categories},
i.e., model categories~$\C$ 
which are enriched over the category $Sp^\Sigma$
of symmetric spectra of~\cite{HSS}, 
compatibly with the smash product and the stable model structure
(see~\cite[Def.~3.5.1]{ss-modules} for details). 
The derived smash product between
symmetric spectra and $\C$ descends to an associative and unital pairing
$$ \sm \ : \ \text{Ho}(Sp^\Sigma) \times \text{Ho }\C \ 
\to \ \text{Ho }\C \ , $$
which can be restricted to the full subcategory $\SW$
of compact objects in Ho$(Sp^\Sigma)$.
A general stable model category is Quillen-equivalent 
to a spectral model category, under mild technical hypothesis
(see~\cite[Thm.~3.8.2]{ss-modules}, 
\cite[Prop.~ 5.5~(a) and 5.6~(a)]{dugger-spectral}
and \cite[Thm.~9.1 and 8.11]{hovey-spectra}
for different sets of sufficient conditions).

In algebraic examples the action of the Spanier-Whitehead category $\SW$ 
and the bounded derived category  $\Dc^b(\mZ)$ (which we used
in the proof of Proposition~\ref{prop-mod n in algebraic})
are related as follows.
The chain functor $C_*:\SW\to D^b(\mZ)$
associates to an object $(X,n)$ of the Spanier-Whitehead category
the reduced singular chain complex of the CW-complex $X$, shifted
up $n$ dimensions. This chain functor is strong symmetric monoidal,
i.e., there are associative, unital and commutative isomorphisms
$C_*(X\sm Y)\iso C_*(X)\tensor^L C_*(Y)$ in $D^b(\mZ)$.
If $\Tc$ is an algebraic triangulated category which is also topological,
then the composite
$$  \SW \times \Tc \ \xra{C_*\times\Id} \
\Dc^b(\mZ) \times \Tc \ \xra{\ \tensor^L } \ \Tc$$
is naturally isomorphic to the smash product pairing.

Now we can make precise in which sense the possibility
of having $n\cdot X/n\ne 0$ in topological triangulated categories
is a 2-local phenomenon.

\begin{prop}\label{prop-away from 2}
If $\Tc$ is a topological triangulated category and $p$ an odd prime, 
then \mbox{$p\cdot X/p=0$} for every object $X$ of $\Tc$.  
 \end{prop}
\begin{proof} The key point is that if $n\not\equiv 2\mod 4$,
then the mod-$n$ Moore spectrum $S/n$ {\em is} annihilated by~$n$
(which is not the case for $n=2$). We recall the easy proof for odd $n$,
which includes all odd primes.
If we take homotopy groups of the defining 
triangle~\eqref{eg-triangle defining Moore} (i.e., apply $[S,-]$, where
$[-,-]$ denotes morphisms in $\SW$) 
we obtain an exact sequence
$$ \pi_1^s \xra{\ n\cdot\ } \pi_1^s \to \pi_1 (S/n) \to 
\pi_0^s \xra{\ n\cdot\ } \pi_0^s \to \pi_0 (S/n) \to 0 \ , $$
using that the stable homotopy groups of spheres vanish in negative dimensions.
Since $\pi_0^s\iso\mZ$ and $\pi_1^s\iso \mZ/2$ we deduce that
$\pi_0(S/n)$ is cyclic of order~$n$ and $\pi_1^s(S/n)=0$ for odd~$n$.
Now we apply $[-,S/n]$ 
to the triangle~\eqref{eg-triangle defining Moore} and
we obtain a short exact sequence
$$ 0 \to [S[1],S/n]\tensor \mZ/n \to 
[S/n,S/n] \to \ \isotope[][n]{[S,S/n]} \to 0 $$
where $\isotope[][n]{A}=\{a\in A\ | \ na=0\}$ 
denotes the group of $n$-torsion points in an abelian group $A$.
By the above calculations, the group $[S[1],S/n]\tensor \mZ/n$ vanishes
and so the group $[S/n,S/n]$ is cyclic of order~$n$ for $n$ odd,
hence $n\cdot S/n=0$.
(For $n=2$ we have $\pi_1(S/2)\iso\mZ/2$ and the analogous
short exact sequence does not split by
Proposition~\ref{prop-mod 2 Moore}. So $[S/2,S/2]\iso\mZ/4$.)

The rest of the argument is then
the same as in Proposition~\ref{prop-mod n in algebraic}.
For an odd prime $p$ we can smash the distinguished triangle
$$ S\ \xra{\ p\cdot\, }\ S\ \xra{\quad}\ S/p\ \xra{\quad}\ S[1]$$
in $\SW$ with the object $X$ and obtain a distinguished triangle in $\Tc$
which shows that $S/p\sm X$ is isomorphic to $X/p$.
Since $p\cdot S/p=0$ in $\SW$ and $\sm$ is biadditive, 
we conclude that $p\cdot X/p= p\cdot(S/p\sm X)=0$.
\end{proof}

We have seen in Proposition~\ref{prop-mod 2 Moore} 
that $p\cdot X/p$ can be non-zero for $p=2$.
On the other hand,
all triangulated categories that I know have the property that 
$p\cdot X/p=0$ for odd primes $p$ and all objects $X$.
 This leaves us with the

\begin{problem}\label{problem geq 3}
Let $p$ be an odd prime. Find a triangulated category $\Tc$ and an object $X$
of $\Tc$ such that $p\cdot X/p\ne 0$, or prove
that in every  triangulated category $\Tc$ we always have $p\cdot X/p= 0$.
\end{problem}

From what we have discussed so far, it is still conceivable that
every topological triangulated category in which~2 is invertible is
algebraic. In particular one can wonder whether the Spanier-Whitehead
category is algebraic after localization at an odd prime. 
We now describe a property of triangulated categories which distinguishes
topological from algebraic examples away from the prime~2.

As before we denote by $K/n$ any cone of $n\cdot\Id_K$,
which comes as part of a distinguished triangle
$$ K\ \xra{\ n\cdot}\ K\ \xra{\ \pi\ }\ K/n\ \xra{\quad}\ K[1] \ . $$
An {\em extension} of a morphism $f:K\to Y$ is then a morphism
$\bar f:K/n\to Y$ satisfying $\bar f\pi=f$.
Such an extension exists if and only if $n\cdot f=0$, and
then the extension will usually not be unique.

\begin{prop}\label{prop-extension in algebraic}
Let $\Tc$ be an algebraic triangulated category, $X$ an object of $\Tc$ 
and $n\geq 2$. Then every morphism $f:K\to X/n$
has an extension $\bar f:K/n\to X/n$  such that some (hence any) mapping cone 
of $\bar f$ is annihilated by $n$.
\end{prop}
\begin{proof}
As in Proposition~\ref{prop-mod n in algebraic}, 
a choice of model for $\Tc$ as the stable category of a Frobenius category
gives a biexact, associative and unital
pairing $\tensor^L:\Dc^b(\mZ)\times\Tc\to\Tc$ and we can take
$X/n=\mZ/n\tensor^L X$ and $K/n=\mZ/n\tensor^L K$. 
We define the extension $\bar f$ as the composite
$$ \mZ/n\tensor^L K\ \xra{\mZ/n\tensor f}\ \mZ/n\tensor^L\mZ/n\tensor^L X
\ \xra{\mu\tensor X} \ \mZ/n\tensor^L X $$
where $\mu:\mZ/n\tensor^L\mZ/n\to\mZ/n$ is the multiplication map
which makes $\mZ/n$ into a ring. 
We choose a distinguished triangle
\begin{equation}\label{eq-triangle}
\mZ/n\tensor^L K\  \xra{\ \bar f\ }\ \mZ/n\tensor^L X 
\ \xra{\ \varphi\ }\ C(\bar f)\ \xra{\ \delta\ }\ \mZ/n\tensor^L K[1]  
\end{equation}
and show that the mapping cone $C(\bar f)$ of $\bar f$ is annihilated by~$n$.

We consider the diagram 
$$
\xymatrix{
\mZ/n\tensor^L\mZ/n\tensor^L K \ar[r]^{\mZ/n\tensor\bar f}
\ar[d]_{\mu\tensor K}&
\mZ/n\tensor^L\mZ/n\tensor^L X \ar[r]^-{\mZ/n\tensor\varphi}
\ar[d]^{\mu\tensor X}& 
\mZ/n\tensor^L C(\bar f) \ar[r]^-{\mZ/n\tensor\delta} \ar@{..>}[d]^\sigma&
\mZ/n\tensor^L\mZ/n\tensor^L K[1]\ar[d]^{\mu\tensor K[1]}\\
\mZ/n\tensor^L K \ar[r]_-{\bar f} &\
\mZ/n\tensor^L X \ar[r]_-{\varphi} & C(\bar f) \ar[r]_-\delta &
\mZ/n\tensor^L K[1]}$$
whose lower row is the distinguished triangle~\eqref{eq-triangle}
and whose upper row is~\eqref{eq-triangle} tensored from the left with $\mZ/n$.
The left square commutes since the multiplication morphism $\mu$ 
is associative in $\Dc^b(\mZ)$.
Since both rows are distinguished triangles, there exists a morphism
$\sigma:\mZ/n\tensor^L C(\bar f)\to C(\bar f)$ making the middle
and right square commute.

We consider the morphism
$$ \pi\tensor C(\bar f)\ :\ C(\bar f)\iso\mZ\tensor^L C(\bar f)
\ \to\  \mZ/n\tensor^L C(\bar f) $$
which satisfies the two relations
\begin{equation}\label{first}
\sigma(\pi\tensor C(\bar f))\varphi  \ =\ 
\sigma(\mZ/n\tensor\varphi)(\pi\tensor\mZ/n\tensor X) \ = \ 
\varphi(\mu\tensor X)(\pi\tensor\mZ/n\tensor X) \ = \ \varphi\\
\end{equation}
and 
\begin{align}
\label{second} \delta\sigma(\pi\tensor C(\bar f))\  &=\ 
(\mu\tensor K[1])(\mZ/n\tensor\delta)(\pi\tensor C(\bar f))\\  
&=\  (\mu\tensor K[1])(\pi\tensor\mZ/n\tensor K[1])\delta \  =\  \delta \ .
\nonumber\end{align}

We claim that the morphism
$$\sigma'\ = \ 2\sigma - \sigma(\pi\tensor C(\bar f))\sigma \ :\
\mZ/n\tensor^L C(\bar f)\to C(\bar f)$$
is a retraction to $\pi\tensor C(\bar f)$. 
Indeed, by~\eqref{first} the composite
of $\varphi:\mZ/n\tensor^L X\to C(\bar f)$ with the morphism 
$\sigma(\pi\tensor C(\bar f))-\Id:C(\bar f)\to C(\bar f)$
becomes trivial, so there exists a morphism $g:\mZ/n\tensor^L K[1]\to C(\bar f)$
such that $g\delta=\sigma(\pi\tensor C(\bar f))-\Id$. But then we have
\begin{align*}
\sigma'(\pi\tensor C(\bar f)) \ &=\ 
2\sigma(\pi\tensor C(\bar f))  - (\Id+g\delta)\sigma(\pi\tensor C(\bar f)) \\
&=\ \sigma(\pi\tensor C(\bar f))  - g\delta\sigma(\pi\tensor C(\bar f)) \
=\ \sigma(\pi\tensor C(\bar f))  - g\delta \ = \ \Id\ , 
\end{align*}
as claimed, where the third equality uses~\eqref{second}.

Now we have shown that $C(\bar f)$ is a direct summand of
$\mZ/n\tensor^L C(\bar f)$, which is annihilated by $n$, and so
we have $n\cdot C(\bar f)=0$.
\end{proof}

In Theorem~\ref{thm-general algebraic} below we prove a generalization
of Proposition~\ref{prop-extension in algebraic}, by a different method.\\

Here is an example showing that in the situation 
of Proposition~\ref{prop-extension in algebraic}
in general there may not exist any extension 
$\bar f:K/n\to X/n$ of $f$ whose mapping cone is annihilated by $n$.
Proposition~\ref{prop-extension in algebraic} and the following 
proposition show that the Spanier-Whitehead category localized at~3
is not algebraic.

We let $\beta_1\in\pi_{10}^s\iso \mZ/6$ be an element of order~3
(so $\beta_1$ generates the 3-primary component of
the 10-dimensional stable stem).
One way to define $\beta_1$ is as the unique element of the Toda bracket
$\td{\alpha_1,\alpha_1,\alpha_1}$ where $\alpha_1\in\pi_3^s\iso\mZ/24$
is the 3-primary part of the element represented by the second Hopf map 
$\nu:S^7\to S^4$.
We let $\tilde\beta_1:S[11]\to S/3$ be any lift in the 
Spanier-Whitehead category $\SW$ of $\beta_1$ to the mod-3 Moore 
spectrum $S/3$. So  $\tilde\beta_1$ is a morphism
whose composite with the connecting map $S/3\to S[1]$ equals the
shift of $\beta_1$. 
We will need below that any such lift $\tilde\beta_1$ is {\em detected 
by $P^3$} in the sense that in any mapping cone
of $\tilde\beta_1$ the Steenrod operation $P^3$ is a non-trivial isomorphism
from mod-3 cohomology in dimension~0 to dimension~12
(see page~60 of~\cite[\S 5]{toda-realizing}).

\begin{prop}\label{prop-mod 3 Moore}
There is no extension of $\tilde\beta_1$ to a morphism 
$\bar\beta:S/3[11]\to S/3$ whose mapping cone is annihilated by~3.
\end{prop}
\begin{proof}
We argue by contradiction and suppose that there exists an extension
$$ \bar\beta\ :\  S/3[11]\ \to\  S/3 $$
of $\tilde\beta_1$ and a distinguished triangle
$$  S/3[11]\ \xra{\ \bar\beta} \  S/3 \ \to \ C(\bar\beta) \ \xra{\ \delta}\
S/3[12] $$
with $3\cdot C(\bar\beta)=0$.

Since the stable stems in dimension~21 and~22 consist only of
torsion prime to~3 we have $\pi_{22}(S/3)=[S[22],S/3]=0$.
So there is a morphism $a:S[23]\to C(\bar\beta)$ lifting $\tilde\beta_1$
in the sense that $\delta a=\tilde\beta_1[12]$.
Since we assumed $3\cdot C(\bar\beta)=0$, the morphism $a$
can be extended to a morphism $\bar a:S/3[23]\to C(\bar\beta)$.
We let $C(\bar a)$ be a mapping cone of $\bar a$ arising as part
of a distinguished triangle
$$  S/3[23]\ \xra{\ \bar a\ } \  C(\bar\beta) \ \to \ C(\bar a\ ) 
\ \xra{\ \delta\ }\ S/3[24]\ . $$
Since the stable stems in dimension~21, 22, 33 and~34 consist only of
torsion prime to~3 we have $\pi_{34}C(\bar\beta)=[S[34],C(\bar\beta)]=0$
and so there is a morphism $b:S[35]\to C(\bar a)$ lifting $\tilde\beta_1$
in the sense that $\delta b=\tilde\beta_1[24]$.

Now we bring cohomology operations into the game to reach a contradiction.
The Moore spectrum $S/3$ has its mod-3 cohomology concentrated in dimensions
0 and 1, where it is 1-dimensional. The mapping cone of $b$ is built
by distinguished triangles from three shifted copies of $S/3$ and one
shifted copy of $S$,
so its mod-3 cohomology is concentrated in dimensions 0, 1, 12, 13, 24, 25
and 36, where it is 1-dimensional.

The morphism $\tilde\beta_1$ is detected by the Steenrod operation $P^3$,
i.e., in the mod-3 cohomology of the mapping cone of $\tilde\beta_1$ 
the operation $P^3$ is non-trivial from dimension~0 to dimension~12.
The stable `cells' (i.e., shifted copies of $S$)
of the mapping cone of $b$ in dimension 12, 24 and 36,
are attached to the two cells directly below by $\tilde\beta_1$, 
so in the cohomology of the mapping cone of $b$,
the 3-fold iterate of $P^3$ is non-trivial from dimension~0 to~36. 
By the Adem relations we have $(P^3)^3=(P^7P^1-P^8)P^1$. Since
the cohomology is trivial in dimension~4 the operation $P^1$
acts trivially, hence so does $(P^3)^3$. 
We have obtained a contradiction, and so
no extension of $\tilde\beta_1$ has a mapping cone which is annihilated by~3.
\end{proof}

In topological triangulated categories, the phenomenon that a morphism
$f:K\to X/n$ may not have any extension whose cone is annihilated by~$n$
is entirely 2- and 3-local, in the following sense. 

\begin{prop}\label{prop-away from 3}
If $\Tc$ is a topological triangulated category and $p$ is a prime 
bigger than~3, then every morphism $f:K\to X/p$ has an extension 
$\bar f:K/p\to X/p$ whose mapping cone is annihilated by $p$.
\end{prop}
\begin{proof} The key point is that if $n$ is prime to~6,
then the mod-$n$ Moore spectrum has an associative multiplication in
the Spanier-Whitehead category
(the multiplication is also commutative, but that 
is not relevant for the current proof).
Again I cannot refrain from giving the simple proof below. 
In contrast, the mod-2 Moore spectrum
does not have a multiplication; the mod-3 Moore spectrum
has a commutative multiplication, but that is not associative in $\SW$, 
see~\cite{toda-realizing}.

As in the proof of
Proposition~\ref{prop-away from 2}, the condition that $n$ is odd 
guarantees that $n\cdot S/n=0$. So there exists a morphism
$\mu:S/n\sm S/n\to S/n$ which splits the two `inclusions'.
We show that $\mu$ is associative. In fact, the  {\em associator}
$$ \mu(\mu\sm\Id-\Id\sm\mu)\ : \ S/n\sm S/n\sm S/n \ \xra{\quad} \ S/n $$
factors as a composite
$$ S/n\sm S/n\sm S/n \ \xra{\delta\sm \delta\sm \delta}\ S[3] \ \to \ S/n  $$
where $\delta:S/n\to S[1]$ is the connecting morphism.
We have $\pi_3^s\iso\mZ/24$ and $\pi_2^s\iso\mZ/2$, so if $n$ is prime to~6 
we have $[S[3],S/n]=\pi_3(S/n)=0$. Thus the associator is trivial or, 
equivalently, $\mu$ is associative in $\SW$.

The rest of the  argument is now essentially the same as in 
Proposition~\ref{prop-extension in algebraic}.
We can smash the distinguished triangle
$$ S\ \xra{\ n\cdot\,}\ S\ \xra{\quad}\ S/n\ \xra{\ \delta\ }\ S[1]$$
in $\SW$ with the object $X$ and obtain a distinguished triangle in $\Tc$
which shows that $S/n\sm X$ is isomorphic to $X/n$.
We use the multiplication $\mu:S/n\sm S/n\to S/n$ to define the extension
$\bar f$ as the composite
$$ S/n\sm K\ \xra{\Id\sm f}\ S/n\sm S/n\sm X\
\xra{\mu\sm\Id}\ S/n\sm X\ . $$
Then we use the same reasoning as in the proof of
Proposition~\ref{prop-extension in algebraic}
to obtain the retraction $\sigma':S/n\sm C(\bar f)\to C(\bar f)$
which shows that $n\cdot C(\bar f)=0$.
\end{proof}

This raises the following question:

\begin{problem}\label{problem geq 5}
Consider a prime $p\geq 5$. Does there exist a triangulated category $\Tc$,
an object $X$ of $\Tc$ and a morphism $f:K\to X/p$ which does not
admit any extension $\bar f$ to $K/p$ whose mapping cone
is annihilated by~$p$?
\end{problem}

\begin{rk}
The proof of Proposition~\ref{prop-away from 3} shows that 
the special features of topological 
over algebraic triangulated categories are closely related to
existence and properties of multiplications on mod-$n$ Moore spectra.
For primes~$p\geq 5$, the mod-$p$ Moore spectrum has a multiplication
in the Spanier-Whitehead category which is commutative and associative.
So on the level of tensor triangulated categories, there does not
seem to be any qualitative difference between the Moore spectrum $S/p$ 
as an object of the Spanier-Whitehead category $\SW$
and $\mZ/p$ as an object of $\Dc^b(\mZ)$, as long as $p\geq 5$.
However, mod-$n$ Moore spectra are never {\em $A_\infty$ ring spectra}, 
but rigorously defining what that means and proving 
it would lead us too far afield. Theorem~\ref{thm-general topological}
below explains how the higher order non-associativity eventually
manifests itself in the triangulated structure of the
Spanier-Whitehead category (i.e., without any reference to the smash product).
\end{rk}

From what we have discussed so far, it is still conceivable that for 
primes $p\geq 5$ the $p$-local Spanier-Whitehead category is algebraic.
We will now introduce an invariant which we then use to
show that this is not the case for any prime $p$.

\begin{defn}
 Consider a triangulated category $\Tc$ and a natural number $n\geq 2$.
 We define the {\em $n$-order} for objects $Y$ of $\Tc$ inductively.
\begin{itemize}
\item Every object has $n$-order greater or equal to~0.
\item For $k\geq 1$, an object $Y$ has $n$-order greater or equal to $k$
if and only if for every object $K$ of $\Tc$ and every morphism $f:K\to Y$ 
there exists an extension $\bar f:K/n\to Y$ such that some (hence any) 
mapping cone of $\bar f$ has $n$-order greater or equal to $k-1$.
\end{itemize}
\end{defn}

We write $n\ord(Y)$, or $n\ord^{\Tc}(Y)$ if we need to
specify the ambient triangulated category, 
for the $n$-order of $Y$, i.e., the largest $k$
(possibly infinite) such that $Y$ has $n$-order greater or equal to $k$.
We define the {\em $n$-order} of the triangulated category $\Tc$
as the $n$-order of some (hence any) zero object,
and denote it by $n\ord(\Tc)$.
We make some observations which are direct consequences of the definitions.

\begin{itemize}
\item The $n$-order for objects is invariant under isomorphism 
and shift.
\item An object $Y$ has positive $n$-order if and only if every morphism
$f:K\to Y$ has an extension to $K/n$, which is equivalent to
$n\cdot f=0$. So $n\ord(Y)\geq 1$ is equivalent to the 
condition $n\cdot Y=0$.
\item  
The $n$-order of a triangulated category is one
larger than the minimum of the $n$-orders of all objects of the form~$K/n$. 
\item Let $\Sc\subseteq\Tc$ be a full triangulated subcategory and
$Y$ an object of $\Sc$. Induction on $k$ shows that if 
$n\ord^{\Tc}(Y)\geq k$, then $n\ord^{\Sc}(Y)\geq k$.
Thus we have 
$n\ord^{\Sc}(Y)\geq n\ord^{\Tc}(Y)$.
In the special case of a zero object we get $n\ord(\Sc)\geq n\ord(\Tc)$.
\item Suppose that $\Tc$ is a $\mZ[1/n]$-linear triangulated category,
i.e., multiplication by $n$ is an isomorphism for every object of $\Tc$.
Then $K/n$ is trivial for every object~$K$ and thus
$\Tc$ has infinite $n$-order. If on the other hand $Y$ is non-trivial, then
$n\ord(Y)=0$. 
\item 
If every object of $\Tc$ has positive $n$-order, then $n\cdot Y=0$
for all objects $Y$ and so $\Tc$ is a $\mZ/n$-linear triangulated category.
Suppose conversely that $\Tc$ is a $\mZ/n$-linear triangulated category.
Then induction on $k$ shows that $n\ord(Y)\geq k$ for all objects
$Y$, and thus every object has infinite $n$-order.
\end{itemize}

The last two items show that the $n$-order is useless if $\Tc$ 
is a {\em $k$-linear} triangulated category for some field~$k$,
since then every $n\in\mZ$ is either zero or a unit in~$k$.

The results which we have obtained so far can be rephrased
using the notion of $n$-order:
if $\Tc$ is algebraic, then Propositions~\ref{prop-mod n in algebraic} 
and~\ref{prop-extension in algebraic} show that 
for every object $X$, the object $X/n$ always has $n$-order
at least~2; we will improve this in Theorem~\ref{thm-general algebraic} below.
Propositions~\ref{prop-mod 2 Moore} and~\ref{prop-mod 3 Moore} show
that in the Spanier-Whitehead category, the Moore spectrum
$S/2$ has 2-order~0 and $S/3$ has 3-order~1;
we generalize this to mod-$p$ Moore spectra
in Theorem~\ref{thm-general topological} below.
If $\Tc$ is topological then for every object $X$, 
the object $X/3$ has $3$-order at least~1, 
by Proposition~\ref{prop-away from 2}.
If $\Tc$ is topological and $p$ is a prime
$\geq 5$, then for every object $X$, the object $X/p$ has $p$-order
at least~2, by Proposition~\ref{prop-away from 3}.

The following theorem generalizes Propositions~\ref{prop-mod n in algebraic}
and~\ref{prop-extension in algebraic}.

\begin{theorem}\label{thm-general algebraic}
  Let $\Tc$ be an algebraic triangulated category and $X$ an object
of $\Tc$. Then for any $n\geq 2$, the object $X/n$ has infinite $n$-order.
In particular, every algebraic triangulated category $\Tc$ has
infinite $n$-order.
\end{theorem}
\begin{proof}(Sketch)
The assumption that $\Tc$ is algebraic gives another piece of
extra structure: for any integer $n$
there exists a triangulated category $\Tc/n$ and
an adjoint pair of exact functors $\rho_*:\Tc\to\Tc/n$
and $\rho^*:\Tc/n\to\Tc$ such that
for every object $X$ of $\Tc$ there exists a distinguished triangle
$$ X\ \xra{\ n\cdot\,}\ X\ \xra{\ \eta\ }\ \rho^*(\rho_* X)\ 
\xra{\quad}\ X[1] $$
where $\eta$ is the unit of the adjunction.

We do not construct $\Tc/n$ in general here, but content ourselves with an
example which gives the main idea. 
If $\Tc=\Dc(A)$ is the derived category of a differential
graded ring~$A$, then we can take $\Tc/n$ as the derived category of the
differential graded ring $A\tensor \overline{\mZ/n}$, 
where $\overline{\mZ/n}$
is a flat resolution of the ring $\mZ/n$, for example the exterior algebra
over $\mZ$ on a generator $x$ of dimension~1 with differential $dx=n$.
The adjoint functor pair $(\rho_*,\rho^*)$ is derived from 
restriction and extension of scalars along the morphism 
$A\to A\tensor \overline{\mZ/n}$.
 
Now we exploit the extra structure to prove the theorem.
Since $X/n$ is isomorphic to $\rho^*(\rho_*X)$ it is enough to show that
for every object $Z$ of $\Tc/n$ and all $k\geq 0$ the object
$\rho^*Z$ has $n$-order greater or equal to~$k$.

We proceed by induction on~$k$; for $k=0$ there is nothing to prove. 
Suppose we have already shown that every $\rho^*Z$ has $n$-order 
greater or equal to~$k-1$ for some positive~$k$. 
Given a morphism $f:K\to \rho^*Z$ in $\Tc$ we can consider
its adjoint $\hat f:\rho_*K\to Z$ in $\Tc/n$; if we apply 
$\rho^*$ we obtain an extension 
$\rho^*(\hat f):K/n\iso \rho^*(\rho_*K)\to\rho^*Z$ of $f$.
We choose a cone of $\hat f$, i.e., a distinguished triangle
$$ \rho_*K \ \xra{\ \hat f\ } \ Z \ \to \ C(\hat f)\ \to\ \rho_*K[1]$$
in $\Tc/n$. Since $\rho^*$ is exact, $\rho^*C(\hat f)$ is a cone
of the extension $\rho^*(\hat f)$ in $\Tc$. By induction, $\rho^*C(\hat f)$ 
has $n$-order greater or equal to $k-1$, which proves that $\rho^*Z$
has $n$-order greater or equal to $k$.
\end{proof}

The next theorem generalizes Propositions~\ref{prop-mod 2 Moore} 
and~\ref{prop-mod 3 Moore} and
shows that topological triangulated categories behave 
quite differently from algebraic ones.

\begin{theorem}\label{thm-general topological}
Let $p$ be a prime. Then in the Spanier-Whitehead category, the mod-$p$ Moore
spectrum $S/p$ has $p$-order  $p-2$. 
Moreover, the  Spanier-Whitehead category, has $p$-order $p-1$.
\end{theorem}

The proof of Theorem~\ref{thm-general topological} has two parts. 
One ingredient is a general statement about
$p$-orders in {\em topological} triangulated categories $\Tc$
which generalizes Propositions~\ref{prop-away from 2}  
and~\ref{prop-away from 3}:
for any object $X$ and prime $p$, the object $X/p$ has $p$-order 
greater or equal to $p-2$. 
The proof of this result
uses the concept and  properties of a {\em coherent action} 
of a mod-$p$ Moore space on an object of a model category, 
see Section~2 of~\cite{schwede-rigid}.
It follows that  any topological triangulated category
(such as the Spanier-Whitehead category)  has $p$-order at least $p-1$.

The second ingredient of Theorem~\ref{thm-general topological}
is the proof that the mod-$p$ Moore spectrum $S/p$ has $p$-order 
at most $p-2$. This uses mod-$p$ cohomology operations and
serious calculational input; in particular,
the proof depends on vanishing results, due to other people,
about the $p$-primary components
of the stable stems in specific dimensions. 
Proposition~\ref{prop-mod 3 Moore} gives the flavor of the proof
which, however, for general primes~$p$ is more involved.
I plan to give a detailed proof elsewhere.

\bigskip

As we just mentioned, every topological triangulated category
has $p$-order at least $p-1$. This leave us with the following question, 
which generalizes Problems~\ref{problem geq 3} and~\ref{problem geq 5}

\begin{problem}
Let $p$ be a an odd prime. Does there exist a triangulated category 
whose $p$-order is strictly less than $p-1$?
\end{problem}

More generally we can ask which values the $n$-orders of
triangulated categories can take.

Now that we have discussed torsion phenomena which can 
distinguish algebraic from topological triangulated categories, 
it is natural to ask whether there are any differences between
topological and algebraic triangulated categories if all primes
are invertible, i.e., in $\mQ$-linear triangulated categories.
The $n$-order is rationally a useless invariant since 
$\mQ$-linear triangulated categories have infinite
$n$-order for all~$n$. 
Similarly, the smash product pairing
of a topological triangulated category with $\SW$ gives no
extra information for $\mQ$-linear triangulated categories since
the chain functor $C_*:\SW\to\Dc^b(\mZ)$ becomes
an equivalence of categories when rationalized
(both sides are in fact rationally equivalent
to the category of finite dimensional graded $\mQ$-vector spaces).

It turns out that rationally the notions of algebraic and
topological triangulated categories essentially coincide. 
Under some mild technical assumptions and cardinality restriction,
every $\mQ$-linear topological triangulated category is algebraic.
More precisely, a theorem of Shipley~\cite[Cor.~2.16]{shipley-DGAs and HZ}
says that every $\mQ$-linear {\em spectral} model category 
(a special kind of stable model category which is enriched over
the stable model category of symmetric spectra) with a set
of compact generators is Quillen-equivalent to dg-modules over
a certain differential graded $\mQ$-category. Thus every triangulated category
equivalent to the homotopy category of a stable model category 
of this kind is algebraic.
I think that the assumption `spectral' is merely of a technical nature
and will eventually be removed. Similarly, I expect that the assumption
of a `set of compact generators' can be relaxed to `well generated'
at the price of allowing localizations of module categories over a
differential graded $\mQ$-category, along the lines of the 
papers~\cite{porta} and~\cite{heider}.
At present, I do not know of a $\mQ$-linear triangulated category which
is not topological.

{\bf Acknowledgement:} 
I would like to thank Andreas Heider for various helpful comments
on this paper.


\begin{thebibliography}{EKMM}

\bibitem[Ad]{adams-generalized}
J.~F.~Adams, {\em Stable homotopy and generalised homology.}
Chicago Lectures in Mathematics. University of Chicago Press, Chicago,
Ill.-London, 1974. x+373 pp.

\bibitem[BF]{BF}
A.~K. Bousfield, E.~M. Friedlander, {\em Homotopy theory of 
{$\Gamma$}-spaces, spectra, and bisimplicial sets.} Geometric applications of 
homotopy theory (Proc. Conf., Evanston, Ill., 1977), II  
Lecture Notes in Math., vol. 658, Springer, Berlin, 
1978, pp.~80--130.


\bibitem[Du]{dugger-spectral}
D.~Dugger, {\em Spectral enrichments of model categories.}
Homology, Homotopy and Applications {\bf 8} (2006), 1--30.

\bibitem[EKMM]{EKMM}
A.~D.~Elmendorf, I.~Kriz, M.~A.~Mandell, J.~P.~May,
{\em Rings, modules, and algebras in stable homotopy theory.
{W}ith an appendix by M.~Cole},
Mathematical Surveys and Monographs, {\bf 47}, American Mathematical Society,
Providence, RI, 1997, xii+249 pp.

\bibitem[He]{heider}
A.~Heider, {\em Two results from Morita theory of stable model categories.}
{\tt arXiv:0707.0707}

\bibitem[Ho1]{hovey-book}
M.~Hovey, {\em Model categories.} Mathematical Surveys and Monographs,
vol.~63, American Mathematical Society,
Providence, RI, 1999, xii+209 pp.

\bibitem[Ho2]{hovey-sheaves}
M.~Hovey, {\em Model category structures on chain complexes of sheaves.}
Trans.\ Amer.\ Math.\ Soc.\ {\bf 353} (2001), 2441--2457.

\bibitem[Ho3]{hovey-spectra}
M.~Hovey, {\em Spectra and symmetric spectra in general model categories.}
J. Pure Appl. Alg. {\bf 165} (2001), 63-127.

\bibitem[HSS]{HSS}
M.~Hovey, B.~Shipley, J.~Smith, {\em Symmetric spectra.}
J.\ Amer.\ Math.\ Soc.\ {\bf 13} (2000), 149--208.

\bibitem[Ke]{keller-differential graded}
B.~Keller, {\em On differential graded categories}. 
Proceedings of the International Congress of Mathematicians, 
Madrid, Spain, 2006, vol.~II, 
European Mathematical Society, 2006, pp.~151--190.

\bibitem[Ly]{lydakis-simplicial}
M.~Lydakis, {\em Simplicial functors and stable homotopy theory.}
{P}reprint (1998).\\ {\tt http://hopf.math.purdue.edu/}

\bibitem[MMSS]{MMSS}
M.~A.~Mandell, J.~P.~May, S.~Schwede, B.~Shipley,
{\em Model categories of diagram spectra.}
Proc.\ London Math.\ Soc. {\bf 82} (2001), 441--512.

\bibitem[Mar]{margolis}
H.~R.~Margolis, {\em Spectra and the Steenrod algebra.
Modules over the Steenrod algebra and the stable homotopy category},
North-Holland Mathematical Library {\bf 29}, North-Holland Publishing Co.,
Amsterdam-New York, 1983, xix+489 pp.

\bibitem[MSS]{nomodel}
F.~Muro, S.~Schwede, N.~Strickland, 
{\em Triangulated categories without models.} 
Invent. Math. {\bf 170} (2007), 231-241.

\bibitem[Ne]{neeman}
A.~Neeman, {\em Triangulated {C}ategories}, Annals of Mathematics Studies,
vol.~148, Princeton University Press, Princeton, NJ, 2001.

\bibitem[Po]{porta}
M.~Porta, {\em The Popescu-Gabriel theorem for triangulated categories.}
{\tt arXiv:0706.4458} 

\bibitem[Q]{Q}
D.~G. Quillen, {\em Homotopical algebra.} Lecture Notes in Math. {\bf 43},
Springer-Verlag, 1967.

\bibitem[Sch]{schwede-rigid}
S.~Schwede, {\em The stable homotopy category is rigid.} 
Annals of Math. {\bf 166} (2007), 837-863.

\bibitem[SS]{ss-modules}
S.~Schwede, B.~Shipley,
{\em  Stable model categories are categories of modules.}
Topology {\bf 42}  (2003),  no. 1, 103--153.

\bibitem[Sh]{shipley-DGAs and HZ}
B.~Shipley, {\em $H\mZ$-algebra spectra are differential graded algebras.}
American J.~Math. {\bf 129} (2007) 351-379. 

\bibitem[SW]{spanier-whitehead}
E.~Spanier, J.~H.~C.~Whitehead, {\em A first approximation to homotopy theory.}
Proc. Nat. Acad. Sci. USA {\bf 39} (1953), 655-660.

\bibitem[To]{toda-realizing}
H.~Toda, {\em On spectra realizing exterior parts of the Steenrod algebra.}  
Topology {\bf 10} (1971), 53--65. 

\bibitem[Vo]{vogt-aarhus}
R.~Vogt, {\em  Boardman's stable homotopy category.}
Lecture Notes Series, No. 21
Matematisk Institut, Aarhus Universitet, Aarhus 1970 i+246 pp.

\end{thebibliography}
\end{document}